\documentclass[12pt,epsfig,tikz,multi]{amsart}
\usepackage{verbatim}
\usepackage{latexsym}
\include{plaquettemacro}
\usepackage{hyperref}
\usepackage{amsrefs}
\usepackage{breqn}
\usepackage{float}
\restylefloat{table}
\usepackage{amsmath,amscd,amssymb}
\usepackage{graphicx}
\usepackage{caption}
\usepackage{subcaption}
\usepackage{enumitem}
\usepackage{xcolor}
\usepackage[utf8]{inputenc}
\usepackage[T1]{fontenc}
 2

\newtheorem{theorem}{Theorem}[section]
\newtheorem{prop}{Proposition}[section]
\newtheorem{definition}[theorem]{Definition}
\newtheorem{example}[theorem]{Example}
\newtheorem{question}[theorem]{Question}
\newtheorem{fact}[theorem]{Fact}
\numberwithin{equation}{section}
\title{Genus Zero Complete Maximal Maps and Maxfaces with an Arbitrary Number of Ends}
\author{Pradip Kumar}
\address{Department of Mathematics, Shiv Nadar Institute of Eminence, Deemed to be University, Dadri 201314, Uttar Pradesh, India.}
\email{pradip.kumar@snu.edu.in}
\author{Sai Rasmi Ranjan Mohanty}
\address{Department of Mathematics, Shiv Nadar Institute of Eminence, Deemed to be University, Dadri 201314, Uttar Pradesh, India.}
\email{sm743@snu.edu.in}
\date{}

\begin{document}

\maketitle
\begin{abstract}
We prove the existence of a genus-zero complete maximal map with a prescribed singularity set and an arbitrary number of simple and complete ends. We also discuss the conditions under which this maximal map can be made into a complete maxface.
\end{abstract}

\section{Introduction}
Similar to minimal surfaces in $\mathbb{R}^3$, maximal surfaces are immersions with zero mean curvature in the Lorentz Minkowski space $\mathbb{E}_1^3$. Maximal surfaces share many similarities with minimal surfaces, such as being critical points of the area functional and admitting the Weierstrass-Enneper representation. However, the theory of complete maximal immersions is limited compared to the theory of complete minimal surfaces. It has been proven \cite{UMEHARA2006} that the plane is the only complete maximal immersion.

Maximal immersions naturally appear with  singularities. Following the terminology from \cite{imaizumi2008}, \cite{LOPEZ20072178}, \cite{UMEHARA2006}, etc., we use the term ``maximal map" to refer to generalized maximal immersions. Maxfaces are a subclass of generalized maximal immersions where the singularities occur only at points where the limiting tangent plane contains a light-like vector \cite{UMEHARA2006}.

The object of the discussion in this article is  genus-zero maximal maps and maxfaces. We summarize the existing results and examples for genus-zero complete maximal maps/maxfaces as follows:
\begin{enumerate}
\item The space-like plane is a maximal immersion.
\item The Lorentzian catenoid is a maxface with two complete, embedded ends.
\item The trinoid is also a genus-zero maxface with three embedded ends \cite{Fujimori2009}.
\item In \cite{Fujimori2009}, the authors provide an example of a genus-zero maxface with two complete ends, with swallowtails, cone-like features, cuspidal edge, and cuspidal crosscaps. These ends are not embedded.
\item Imaizumi and Kato \cite{imaizumi2008} discuss different types of simple ends and classify maximal surfaces of genus zero with at most three embedded ends.
\end{enumerate}

In \cite{imaizumi2008}, authors discuss the classifications of ends of the maximal map, and in \cite{Ando},  we see the regularity results about the ends of zero mean curvature surfaces in the Lorentz-Minkowski space.

The question addressed in this article is as follows:
\begin{question}\label{question}
Given a prescribed singular set and the prescribed nature of singularity, can we construct a \textbf{genus-zero} maximal map or maxface with an \textbf{arbitrary number }of \textbf{simple ends} and \textbf{complete ends}?
\end{question}

One possible approach to constructing a maximal map or maxface is by starting with a minimal immersion. For a minimal immersion, there exists a corresponding maximal immersion (known as a companion \cite{UMEHARA2006}). However, this method does not always yield a complete maxface.  An example of this is the Jorge-Meeks minimal surface \cite{Meeks1983}.

Weierstrass data for Jorge-Meek's minimal surface is the following:
$$
g(z)=z^n\mbox{ and } \omega=\dfrac{dz}{(z^{n+1}-1)^2}
.$$
This gives a complete minimal surface  on $\mathbb C\cup\{\infty\}$ with punctures $\{1, \zeta, \zeta^2,... ,\zeta^{n-1}\}$, where $\zeta= e^{\frac{2i \pi}{n}}$. 
We see $(g_0= i g, \omega_0=\omega)$ gives a maximal map, but this is not complete maxface since the singular set is not compact.

To construct a complete maximal map or maxface, the Bj\"{o}rling problem and its solution are not always sufficient. In some cases, starting from the Bj\"{o}rling data leads to a non-complete maxface. Examples in \cite{Kim2007} (for instance, see Example 3.11) do not provide the Weierstrass data for complete maxfaces. 

To construct examples of complete maxfaces, we need to find a suitable meromorphic function $g$ and a holomorphic 1-form $\omega$ on $\mathbb{C}\cup\{\infty\}\setminus\{p_1,...,p_n\}$ that satisfies the period condition (as in (3) of Theorem \ref{weirstrassdata}). However, unlike in the case of minimal surfaces, when solving the balance equation of the period problem, we have to ensure that the meromorphic function $g$ is chosen in a way that $|g(p_i)|\neq 1$ after extension to the ends. Due to this restriction, most methods and techniques used for minimal surfaces can not be directly applied.

In this article, we prove the existence of maximal maps or maxfaces with complete ends, and since we are not specifically looking for embedded ends, the method for solving the period problem can be simplified. Here, we explain the method we will use.

For any $r$ and the subset $\{p_1, ..., p_r\}\subset \mathbb{C}$, if we take
$$
f(z)=\frac{dz}{(z-p_1)^2\cdots (z-p_r)^2}; \mbox{ and } g(z)=\sum_{i=0}^{r-1}a_i z^i+\sum_{j=1}^{r}\frac{b_j}{z-p_i},
$$
then it is well-known (see \cite{KichoonBook}) that there exist constants $a_i$'s and $b_j$'s such that $(g,f(z)dz)$ gives a complete branched minimal surface on $\mathbb{C}\cup \{\infty\}\setminus\{p_1, ..., p_r\}$. The same data and solution can be used to obtain a maximal map with any number of ends. However, this solution does not guarantee that $|g(p_i)|\neq 1$ at the ends, and thus we do not obtain complete ends. Moreover, the $a_i$'s and $b_j$'s are not explicitly determined, so we do not have control over the exact singularity set for the maxface. In this article, we will modify this method to solve the period problem and find the desired maxface or maximal map.

In Section 3, we explain the period problem and prove the following results:
\begin{enumerate}
\item Let $C$ be a singular curve in $\mathbb{C}$ as defined in the definition  \ref{singularcurve}, and let $m, n \in \mathbb N^*$.  Then:
\begin{enumerate}
\item (Theorem \ref{intro1}) there exists a maximal map with $m$ complete ends and $C$ as its singularity set.
\item (Theorem \ref{lastthm}) there is a maximal map with $m$ complete ends, $n$ simple ends, and $C$ as its singularity set. All simple ends, in this case, are embedded.
\end{enumerate}
\item (Theorem \ref{intro2}) Let $X: \mathbb{C}\cup\{\infty\}\setminus\{p_1,...,p_n\}\to \mathbb{E}_1^3$ be a complete maximal map with the Gauss map $g$. Then, there exists a finite set $E$ and a family of complete maxfaces on $M\setminus E$ with the Gauss map $\tilde{g}$, such that $\tilde{g}$ is very close to $g$ away from the ends.
\end{enumerate}

It is important to note that these results are trivially solved in the case of minimal immersions (such as the Jorge-Meeks surface). However, for maximal surfaces, due to the presence of non-isolated singularities, finding such results requires a slightly different approach.

\section{Preliminary}
We denote the Lorentz Minkowski space as $\mathbb{E}_1^3$, which is the vector space $\mathbb{R}^3$ with the metric $dx^2+dy^2-dz^2$. Similar to the case of minimal immersions, Esutidillo and Romero \cite{Estudillo1992} have proved the Weierstrass-Enneper representation for generalized maximal immersions (maximal maps). We mention it here.

\begin{theorem}[Weierstrass-Enneper representation for maximal maps \cite{Estudillo1992}]\label{weirstrassdata}
Let $M$ be a Riemann surface, $g$ be a meromorphic function, and $\omega= f(z)dz$ be a holomorphic 1-form on $M$ such that:
\begin{enumerate}
\item If $p\in M$ is a pole of $g$ of order $m$, then $f$ has a zero at $p$ of order at least $2m$.
\item $|g|\not\equiv 1$ on $M$.
\item $Re\int_{\gamma}((1+g^2)\omega, i(1-g^2)\omega, -2g\omega)=0$ for all closed loops $\gamma$ on $M$.
\end{enumerate}
Then, the map $X: M\to \mathbb{E}_1^3$ given by
\begin{equation}
X(p)=Re\int_{0}^{p}((1+g^2)\omega, i(1-g^2)\omega, -2g\omega)
\end{equation}
is a maximal map with base point $0\in M$.

Moreover, any maximal map can be expressed in this form.
\end{theorem}

The third condition in the above theorem \ref{weirstrassdata} is called the period condition. For the maximal map $X:M\to \mathbb{E}_1^3$ with the Weierstrass data $(g, \omega)$, the pullback metric is given by $ds^2=(1-|g|^2)^2|\omega|^2$.

The singularities of the maximal map $X: M\to \mathbb{E}_1^3$ are given by the set $Sing(X)= \{p\in M: |g(p)|=1 \text{ or } \omega(p)=0\}$. We call the set $\{p\in M: |g(p)|=1\}$ the singularity set of the maximal map and $\{p\in M: \omega(p)=0\}$ the set of branch points.

For the maxface, there are no branch points, and we have a singularity where $|g|=1$ and $\omega\neq 0$. Therefore, if $(g,\omega)$ is the Weierstrass data for the maxface, we must have
\begin{equation}\label{eqndivisorcond}
(\omega)_{0}=2(g)_{\infty} \text{ on } M.
\end{equation}
We refer to the above equation as the \textbf{divisor condition}.

The completeness of a maxface can be understood using the following equivalent criteria:
\begin{fact}
A maxface is complete (see \cite{UMEHARA2006}, corollary 4.8) if and only if:
\begin{enumerate}
\item $M$ is bi-holomorphic to $\overline{M}\setminus\{p_1,...,p_n\}$.
\item $|g|\neq 1$ at $p_i$, $1\leq i\leq n$.
\item The induced metric $ds^2$ is complete at the ends.
\end{enumerate}
\end{fact}
In this case, the points $p_i$'s are called ends for the maxface, and $g$ and $\omega$ can be extended as meromorphic function and meromorphic 1-form respectively at the ends.

Similarly, for the maximal map defined on $\overline{M}\setminus \{p_1, ..., p_n\}$, we say it is a complete maximal map if and only if $|g(p_i)|\neq 1, 1\leq i\leq n, $ and the induced metric $ds^2$ is complete at the ends.

In the next section, we will prove the existence of a maximal map with an arbitrary number of complete ends and a prescribed singularity set. We will start the discussion by finding a complete maximal map with the prescribed Gauss map.
\section{Genus zero complete maxface with prescribed singularities}
We will start this section by constructing a maximal map with a given Gauss map.

\subsection{Maximal map with an arbitrary number of complete ends and a prescribed Gauss map}\label{section3.1}
Let $g$ be a meromorphic function on $\mathbb{C}\cup\{\infty\}$ such that
\begin{equation}\label{gaussmap1}
(g)_{\infty}=\sum_{i=1}^{n}x_ip_i, \quad x_i\in \mathbb{N}^*,  1\leq i\leq n.
\end{equation}
Without loss of generality, we assume $p_n=\infty$.

Next, we will find a suitable holomorphic 1-form $\omega$ such that $(g,\omega)$ becomes the Weierstrass data for a maximal map $X$ defined on $\mathbb{C}\cup\{\infty\}\setminus\{p_1,...,p_n\}$, and all the ends are complete.

Consider the following expression for $f(z)$:
\begin{equation}\label{heightDifferential1}
f(z)= \sum_{i=1}^{n-1}\frac{a_i}{(z-p_i)^2}+\sum_{j=0}^{2n-2} b_j z^j,
\end{equation}
where $a_i, b_j\in \mathbb{C}$ (we will find suitable values for $a_i$'s and $ b_j$'s). Set $\omega_0:=f(z)dz$.

We need to show that there exist constants $(a_i; b_j)\in \mathbb{C}^{3n-2}$ such that the period condition, as in the theorem \ref{weirstrassdata}, holds for all loops on $\mathbb{C}\cup\{\infty\}\setminus\{p_1,..., p_n\}$. Since $\omega_0$ is an exact form, and if $g^2\omega_0$ and $g\omega_0$ have no residues on $\mathbb{C}\cup\{\infty\}$, the pair $(g, \omega_0)$ satisfies the period condition. Therefore, we need to find values of $(a_i; b_j)$ such that the following equations hold:
\begin{equation}\label{eqnresgw}
\text{Res}_{p_i}(g\omega_0)= 0, \quad \text{Res}_{p_i}(g^2\omega_0)= 0, \quad 1\leq i\leq n.
\end{equation}
Equation \ref{eqnresgw} gives a homogeneous system of linear equations in the variables $a_i$ and $b_j$. Moreover, since the total residues of $g\omega_0$ and $g^2\omega_0$ are zero on $\mathbb{C}\cup\{\infty\}$, we have a homogeneous system of at most $2n-2$ linear equations in $3n-2$ variables $a_i$ and $b_j$. Therefore, the solution space (denoted by $Sol\subset \mathbb{C}^{3n-2}$) for the system of equations \ref{eqnresgw} has a dimension of at least $n$. This proves that we can choose suitable $(a_i; b_j)\in Sol\setminus\{0\}$ such that $(g,\omega_0)$ satisfies the period condition.

Since $\omega_0$ has poles at $p_i$'s, the divisor condition for $(g,\omega_0)$, as stated in the theorem \ref{weirstrassdata}, holds trivially on $\mathbb{C}\cup\{\infty\}\setminus\{p_1,...,p_n\}$. Moreover, since the solution space is of dimension $n$, for $(a_i; b_j)\in Sol\setminus\{0\}$, at least one of the $b_j$'s is non-zero. Therefore, $p_n$ is a pole of both $\omega_0$ and $g$, making $p_n$ a complete end for the maximal map.

However, it may happen that for some $i$ ($1\leq i\leq n-1$), $a_i=0$, in which case $p_i$ may be a zero of $\omega_0$. Thus, while the pair $(g,\omega_0)$ gives a maximal map, $p_i$ may not be a complete end for this maximal map. To address this, we aim to modify $\omega_0$ such that the resulting pair becomes the Weierstrass data for a maximal map with all $p_i$'s as poles of $\omega_0$, ensuring that all $p_i$'s are complete ends.

To achieve this, for $1\leq i\leq n-1$, we consider
$$
f_i=\sum_{k=1}^{2n}\frac{\alpha_k}{(z-p_i)^{2+k}}
$$ and we want to find $(\alpha_k)\in \mathbb C^{2n}$ with at least one $\alpha_k\neq 0$ such that $(g,f_i)$ satisfies the period condition on $\mathbb{C}\cup\{\infty\}$. By using the technique mentioned earlier, i.e., $Res_{p_j}(gf_idz)=0$ and $Res_{p_j}(g^2f_idz)=0$ for $1\leq j\leq n$, we obtain $2n-2$ independent linear homogeneous equations with $2n$ variables. Therefore, there exists $(\alpha_k)\in \mathbb{C}^{2n}\setminus\{0\}$ such that $(g,f_i)$ satisfies the period condition and $p_i$ is a pole of $f_i$.

Let $F:=f+\sum_{i=1}^{n-1}f_i$ and ${\omega}:=F(z)dz$. It is clear that $g$ and ${\omega}$ satisfy the divisor condition and the period condition on $\mathbb{C}\cup\{\infty\}\setminus\{p_1,...,p_n\}$. Moreover, all $p_i$'s are poles of $g$ and ${\omega}$, making all ends complete.

To have a complete maximal map $X: \mathbb{C}\cup\{\infty\}\setminus\{p_1,...,p_n\}\to \mathbb{E}_1^3$ with the Weierstrass data $(g,{\omega})$, we only need to prove that the singular set is compact. It is evident from the fact that $g$ has poles at the ends and the manner in which we adjusted $\omega_0$ to obtain $\omega$, which implies that $\{z: |g|=1 \}\cup\{z: {\omega}(z)=0\}$ is compact.

To summarize, we have proved the following proposition:

\begin{prop}\label{completemaximalmap with prescribed gaussmap}
Let $g$ be a meromorphic function on the Riemann sphere such that the number of distinct poles of $g$ is $n$. Then there exists a maximal map with $n$ complete ends having the Gauss map $g$.
\end{prop}

The maximal map obtained from the proposition \ref{completemaximalmap with prescribed gaussmap} is not necessarily a maxface because the 1-form $\omega$ may have some zeros on $\mathbb{C}\cup\{\infty\}\setminus\{p_1,...,p_n\}$. In subsection \ref{sec:perturbationmaxface}, we will perturb the Gauss map to obtain the complete maxface. Furthermore, the above method takes care of the ends as complete, but in this process, we obtain higher-order poles at the ends, which prevents the possibility of embedded ends.

Moreover, if we replace $F$ as defined above by

\begin{equation}\label{eq:heightDifferential2}
F(z)=\sum_{i=1}^{2n}a_ih^i,
\end{equation}

where, for $n\geq 2$, $h(z)=\prod_{j=1}^{n-1}\frac{(z-c)^{2n-2}}{(z-p_j)^2}$ and $c\neq p_j$ for $1\leq j \leq n$, and for $n=1$, $h(z)$ is any polynomial, we obtain the same result. Additionally, there are other ways to define the required $\omega$. In all these cases, we can increase the number of variables to obtain a larger dimensional solution space, which may help us in obtaining the maxface with a prescribed nature of singularity. This aspect needs to be explored further.  

Below, we will define the singular curve and discuss the possibility of the existence of a maximal map with such a curve as its singularity set. 

\subsection{Maximal map with an arbitrary number of  complete ends and a prescribed singularity set}

For a given curve or collection of disjoint union of curves $C\subset \mathbb C$, we aim to find a complete maximal map with an arbitrary number of complete ends such that $C$ is the singularity set. So firstly, for a given $C$, we will be looking for the meromorphic function $g$ on $\mathbb C\cup\{\infty\}$ such that for $z\in C$, $|g(z)|=1$, and then applying the proposition \ref{completemaximalmap with prescribed gaussmap}, we will get the required maximal map. 

Considering the above, we will take $C$ as follows in the rest of the article. 
\begin{definition}\label{singularcurve}
    The singular curve $C$ is a curve or disjoint union of curves such that there exists a meromorphic function $g$ on $\mathbb{C} \cup \{\infty\}$ such that $g(C)$ is mapped to the circle $|z|=R$. 
\end{definition}

There are many examples of singular curves. Below, we present two examples that constitute a significant class of singular curves.
\begin{example}
We begin by recalling the definition of the Schwarzian derivative of a meromorphic function $f(z)$;
$S[f](z):= \left(\frac{f^{\prime\prime}(z)}{f'(z)}\right)^\prime - \frac{1}{2} \left(\frac{f^{\prime\prime}(z)}{f'(z)}\right)^2.$  The M\"{o}bius transformations are the only functions for which $S[f](z)=0$ \cite{Osgoos}.

Consider a regular simple closed curve $C$ such that the Schwarzian derivative of any conformal parameterization $f(z)$ of $C$ is identically zero. In this case, there exists a meromorphic function $G(z)$ on $\mathbb{C} \cup \{\infty\}$ that maps $C$ to the unit circle.
\end{example}

\begin{example}
    Another example of a singular curve is the subset $C = |G|^{-1}(1)$, where $G$ is defined as follows:

$$G(z)= \frac{(z-1)(z^2+3 z+1)}{(z+1)(z^2-3 z+1)}.$$  This curve and the meromorphic function $G$ are taken from the proof of the theorem A \cite{Fujimori2009}, where the authors investigated various types of singularities on maxfaces. For this particular $G$, the singular curve $C$ consists of three disjoint (topological) circles on the Riemann sphere, including the imaginary axis.
\end{example}

Moreover, any meromorphic function $G$ that maps a circle of radius $R$ to the unit circle is a Blaschke product, as discussed by Gronwall in \cite{Gronwall1912}, and is given by
\begin{equation}\label{gaussmap2}
G(z)=\prod_{i=1}^{m}\frac{R(z-a_i)}{a_i(z-a_i^\prime)}\prod_{j=1}^{n}\frac{b_j(z-b_j^\prime)}{R(z-b_j)},
\end{equation}
where $a_i$, $1\leq i\leq m$, and $b_j$, $1\leq j \leq n$, are points inside the circle $|z|=R$, and $a_ia_i^\prime=R^2$ and $b_jb_j^\prime=R^2$. For $a_i=0$, the factor $\frac{R(z-a_i)}{a_i(z-a_i^\prime)}$ in the equation \ref{gaussmap2} is replaced by $-\frac{z}{R}$, and  for $b_j=0,$ the factor $\frac{b_j(z-b_j^\prime)}{R(z-b_j)}$ is replace by $-\frac{R}{z}$. Note that $a_i$'s and $b_j$'s are repeated according to the order of zeros and poles.

For the given singular curve $C$, as defined in the definition \ref{singularcurve}, there is a meromorphic map $f: \mathbb C\cup\{\infty\}\to \mathbb C\cup\{\infty\}$ such that $f(C)$ is the circle of radius $R$.  We choose $G$ as described in the equation \ref{gaussmap2} and compose it with the meromorphic map $f$. This composite map, denoted as $g$, satisfies $|g(z)|=1$ on $C$. Additionally, by choosing $G$ appropriately, we have the flexibility to increase the number and order of poles of $g$.

Furthermore, after the change of coordinates, we assume that one of the poles of $g$ is at $\infty$. We denote the poles of $g$ as $p_i$, where $1\leq i \leq n$, and write $(g)_{\infty}=\sum_{i=1}^{n}x_i p_i$, as in the equation \ref{gaussmap1}.

Now, if we start with the meromorphic function $g$ as described in the previous paragraph, we can use the proposition \ref{completemaximalmap with prescribed gaussmap} to establish the following theorem:

\begin{theorem}\label{intro1}
Given a singular curve $C$ as in the definition \ref{singularcurve}, and $m\in \mathbb N^*$, there exists a maximal map with $m$ complete ends and $C$ as its  singularity.
\end{theorem}

Theorem \ref{intro1} establishes the existence of the desired maximal map. However, if we want a maxface, we can get it by slightly modifying the normal vector (the Gauss map $g$). Specifically, in the following subsection, given a complete maximal map with Weierstrass data $(g, \omega)$, we will find a suitable meromorphic function $g_0$ such that $(\tilde{g} = g + g_0, \omega)$ yields a complete maxface. Moreover, $\tilde{g}$ can be chosen arbitrarily close to $g$, away from the ends.

\subsection{Perturbation of the Weierstrass data for the complete maxface}\label{sec:perturbationmaxface}

Let $(g, \omega = f dz)$, as defined in the equations \ref{gaussmap1} and \ref{eq:heightDifferential2} respectively, be the Weierstrass data for the maximal map $X: \mathbb{C}\cup\{\infty\}\setminus\{p_1,\ldots,p_n\} \to \mathbb{E}_1^3$ with complete ends, as stated in the theorem \ref{intro1}.

Let $\{z_j: 1\leq j\leq l\}$ be the zeros of $\omega$ and $Ord_{z_j}\omega = m_j > 0$. With $\alpha_i\in\mathbb{C}$, we define
\begin{align}\label{g_0}
g_0(z) = \frac{1}{\prod_{j=1}^{l}(z-z_j)^{2m_j}}\left(\sum_{i=1}^{3(l+n)}\alpha_i z^i\right).
\end{align}

We take $\tilde{g} = g + g_0$. In the following, we will determine $(\alpha_i)\in\mathbb{C}^{3(l+n)}$ such that $(\tilde{g}, \omega)$ serves as the Weierstrass data for the complete maxface on $\mathbb{C}\cup\{\infty\}\setminus\{p_1,\ldots,p_n,z_1,\ldots,z_l\}$, with $\tilde{g}$ being close to $g$ away from the ends.

The divisor condition, as given in the equation \ref{eqndivisorcond}, trivially holds for it to be a maxface.

If the period condition holds for $(\tilde{g}, \omega)$, then the corresponding maxface will be complete. This follows from the construction, since we observe that $|\tilde{g}|\neq 1$ at the ends $\{p_i; z_j : 1\leq i\leq n, 1\leq j\leq l\}$. Additionally, at the end, the metric $ds^2$ is complete. Hence, all ends are complete.

Thus, we need to find $(\alpha_i)$ such that the period condition for $(\tilde{g}, \omega)$ is satisfied. Since $g\omega$ and $g^2\omega$ have no residues on $\mathbb{C}\cup\{\infty\}$ and $\omega$ is exact, we need to find $(\alpha_i)$ such that $g_0\omega$, $g_0^2\omega$, and $gg_0\omega$ have no residues on $\mathbb{C}\cup\{\infty\}$. This can be achieved by employing the same technique as in subsection \ref{section3.1}.

Moreover, for any $\epsilon > 0$, the vector $\epsilon \alpha = (\epsilon \alpha_i)$ is also a solution. Therefore, $g_0$ can be suitably chosen such that away from the ends $z_j$, we have $|g_0| < \epsilon$. This proves that $|\tilde{g} - g| = |g_0|$ is very small away from the ends.

Combining all the above, we establish the following result:

\begin{theorem}\label{intro2}
Let $X: M = \mathbb{C}\cup\{\infty\}\setminus\{p_1,\ldots,p_n\} \to \mathbb{E}_1^3$ be a complete maximal map with the Gauss map $g$. Then, there exists a finite set $E$ and a family of complete maxfaces on $M\setminus E$ with the Gauss map $\tilde{g}$ such that $\tilde{g}$ is very close to $g$ away from the ends.
\end{theorem}

In all of the above proofs, we solved the period problem by increasing the order of poles and variables sufficiently. This ensured the existence of complete ends, but it did not guarantee simple ends. However, with some modifications and a similar technique, it is possible to obtain a maximal map with the prescribed simple ends. We explain the approach here.

\subsection{Maximal map with a given Gauss map and an arbitrary number of simple ends}

Let $g$ be a meromorphic function as defined in the equation \ref{gaussmap1}, and for any $q_j$ (where $1\leq j\leq m$) on $\mathbb{C}$ distinct from $p_i$'s, we aim to find a suitable holomorphic 1-form $\omega$ such that $(g,\omega)$ serves as the Weierstrass data for the complete maximal map on $\mathbb{C}\cup\{\infty\}\setminus\{p_1,\ldots,p_n,q_1,\ldots,q_m\}$, where $p_i$'s are complete ends and $q_j$'s are simple (in particular, embedded) ends.

We introduce some $a_i$'s and, we define
\begin{equation}\label{eqn:f:forsimpleend}
f(z) = \frac{1}{\prod_{j=1}^{m}(z-q_j)^2}\left(\sum_{i=1}^{2m+2n} a_i h^i\right),
\end{equation}
where for $n\geq 2$, $h(z)=\prod_{k=1}^{n-1}\frac{(z-c)^{2n-2}}{(z-p_k)^2}$ and $c\neq p_i$ for $1\leq i \leq n$, and for $n=1$, $h(z)$ is any polynomial.

By employing a similar technique as discussed in the section \ref{section3.1}, we can prove that there exists a choice of $(a_i)$ such that $(g,f(z)dz)$ provides the desired maximal map.

Formally, we state the theorem as follows:

\begin{theorem}\label{lastthm}
Given a meromorphic function $g$ and $m\in\mathbb N^*$, there exists a maximal map with complete ends at the poles of $g$ and $m$ simple ends. All simple ends, in this case, are embedded.
\end{theorem}

Again, if we examine the proof, we realize that the choice of $f(z)$ is flexible. We have the freedom to increase the variables and points to obtain any desired number of ends.

\section{General remarks and conclusion}

The main results of this article are presented in theorems \ref{intro1}, \ref{intro2}, and \ref{lastthm}. These theorems establish the existence of maxfaces/maximal maps with a prescribed singular set and an arbitrary number of complete ends, as well as simple ends. However, finding such non-complete maxfaces is straightforward. For example, if we take $a, b, c, d$ and $e$ as distinct real numbers and set $g(z)=e^{iz}$ and $f(z)=((z-b)(z-d)^2-i(z-c)(z-e)^2)e^{-iz}$, then $(g,f(z)dz)$ forms a Weierstrass data for a non-complete maxface on $\mathbb{C}$. An interesting observation, in this case, is that the singularities at $a, b, c, d$, and $e$ correspond to cuspidal-edge, swallowtails, cuspidal-butterflies, cuspidal $S_1^-$, and cuspidal-crosscaps, respectively.

Similar questions can be posed regarding the existence of complete maxfaces with various prescribed types of singularities occurring at different points. As the solution space for the required constants is vast, it is likely possible to choose suitable constants. However, this requires an improved method for solving the period problem.
\section{Acknowledgments}
The authors would like to express their sincere gratitude to the anonymous reviewers for their valuable comments and suggestions, which greatly contributed to the improvement of this article. 

\medskip

\bibliography{ref}

\begin{bibdiv}
\begin{biblist}

\bib{Estudillo1992}{article}{
      author={Estudillo, F. J.~M.},
      author={Romero, A.},
       title={Generalized maximal surfaces in lorentz–minkowski space l3},
        date={1992},
     journal={Mathematical Proceedings of the Cambridge Philosophical Society},
      volume={111},
      number={3},
       pages={515–524},
}

\bib{Fujimori2009}{article}{
      author={Fujimori, S.},
      author={Rossman, W.},
      author={Umehara, M.},
      author={Yamada, K.},
      author={Yang, S.-D.},
       title={New maximal surfaces in minkowski 3-space with arbitrary genus
  and their cousins in de sitter 3-space},
        date={2009},
        ISSN={1420-9012},
     journal={Results in Mathematics},
      volume={56},
      number={1},
       pages={41},
         url={https://doi.org/10.1007/s00025-009-0443-4},
}

\bib{Gronwall1912}{article}{
      author={Gronwall, T.~H.},
       title={On analytic functions of constant modulus on a given contour},
        date={1912/13},
        ISSN={0003-486X},
     journal={Ann. of Math. (2)},
      volume={14},
      number={1-4},
       pages={72\ndash 80},
         url={https://doi.org/10.2307/1967601},
      review={\MR{1502440}},
}

\bib{imaizumi2008}{article}{
      author={Imaizumi, T.},
      author={Kato, S.},
       title={Flux of simple ends of maximal surfaces in ${ R}^{2,1}$},
        date={2008--},
     journal={Hokkaido Math. J.},
      volume={37},
      number={3},
       pages={561\ndash 610},
         url={https://doi.org/10.14492/hokmj/1253539536},
}

\bib{Meeks1983}{article}{
      author={Jorge, Luqu\'{e}sio~P.},
      author={Meeks, William~H., III},
       title={The topology of complete minimal surfaces of finite total
  {G}aussian curvature},
        date={1983},
        ISSN={0040-9383},
     journal={Topology},
      volume={22},
      number={2},
       pages={203\ndash 221},
         url={https://doi.org/10.1016/0040-9383(83)90032-0},
      review={\MR{683761}},
}

\bib{Kim2007}{article}{
      author={Kim, Y.~W.},
      author={Yang, S.-D.},
       title={Prescribing singularities of maximal surfaces via a singular
  björling representation formula},
        date={2007},
        ISSN={0393-0440},
     journal={Journal of Geometry and Physics},
      volume={57},
      number={11},
       pages={2167 \ndash  2177},
}

\bib{LOPEZ20072178}{article}{
      author={López, Rafael},
       title={On the existence of spacelike constant mean curvature surfaces
  spanning two circular contours in minkowski space},
        date={2007},
        ISSN={0393-0440},
     journal={Journal of Geometry and Physics},
      volume={57},
      number={11},
       pages={2178\ndash 2186},
  url={https://www.sciencedirect.com/science/article/pii/S0393044007000769},
}

\bib{Ando}{article}{
      author={Naoya, Ando},
      author={Kohei, Hamada},
      author={Kaname, Hasimoto},
      author={Shin, Kato},
       title={{Regularity of ends of zero mean curvature surfaces in
  $\mathbf{R}^{2,1}$}},
        date={2022},
     journal={Journal of the Mathematical Society of Japan},
      volume={74},
      number={4},
       pages={1295 \ndash  1334},
         url={https://doi.org/10.2969/jmsj/85018501},
}

\bib{Osgoos}{article}{
      author={Osgood, Brad},
       title={Old and new on the schwarzian derivative},
        date={1998},
       pages={275\ndash 308},
         url={https://doi.org/10.1007/978-1-4612-0605-7_16},
}

\bib{UMEHARA2006}{article}{
      author={Umehara, M.},
      author={Yamada, K.},
       title={Maximal surfaces with singularities in minkowski space},
        date={2006},
     journal={Hokkaido Math. J.},
      volume={35},
      number={1},
       pages={13\ndash 40},
         url={https://doi.org/10.14492/hokmj/1285766302},
}

\bib{KichoonBook}{book}{
      author={Yang, Kichoon},
       title={Complete minimal surfaces of finite total curvature},
      series={Mathematics and its Applications},
   publisher={Kluwer Academic Publishers Group, Dordrecht},
        date={1994},
      volume={294},
        ISBN={0-7923-3012-9},
         url={https://doi.org/10.1007/978-94-017-1104-3},
      review={\MR{1325927}},
}

\end{biblist}
\end{bibdiv}
 \end{document}